\documentclass[12pt,a4paper]{article}
\usepackage{amsmath,amsfonts,amssymb,amsthm,amscd}
\usepackage{hyperref}
\usepackage[T1]{fontenc}
\usepackage{color}
\usepackage{epstopdf}

\usepackage{braket}
\usepackage{xypic}
\usepackage{amsmath}

\usepackage[american]{babel}
\usepackage{color} 
\usepackage{epsfig}
\usepackage{caption}
\usepackage{rotating}
\usepackage{setspace}
\usepackage{fancyhdr}
\usepackage{booktabs}
\usepackage{mathrsfs}
\usepackage{amsthm}
\usepackage{amsmath,amssymb}
\usepackage{enumerate}
\usepackage{tikz}
\usepackage{youngtab}
\usepackage{hyperref}
\usepackage{tensor}
\usepackage{mathabx}
\usepackage{booktabs}
\usepackage{mathrsfs}
\usepackage{subfigure}
\usepackage{amsmath,amssymb}

\newcommand{\ra}{\rightarrow}

\newcommand{\zz}{\mathbb{Z}}
\newcommand{\nn}{\mathbb{N}}

\newcommand{\rr}{\mathbb{R}}

\DeclareMathOperator{\PSL}{PSL}

\DeclareMathOperator{\Aut}{Aut}

\DeclareMathOperator{\supp}{supp}
\DeclareMathOperator{\SL}{SL}
\DeclareMathOperator{\wsr}{wsr}
\DeclareMathOperator{\Ped}{Ped}

\DeclareMathOperator{\asr}{asr}

\newtheorem{thm}{Theorem}[section]
\newtheorem{pro}[thm]{Proposition}
\newtheorem{cor}[thm]{Corollary}
\newtheorem{lem}[thm]{Lemma}

\theoremstyle{definition}
\newtheorem{rem}[thm]{Remark}
\newtheorem{defn}[thm]{Definition}

\begin{document}

\title{On $C^*$-algebras associated to actions of discrete subgroups of $\SL(2,\rr)$ on the punctured plane}

\author{Jacopo Bassi} 

\maketitle

\begin{abstract}
\noindent Dynamical conditions that guarantee stability for discrete transformation group $C^*$-algebras are determined. The results are applied to the case of some discrete subgroups of $\SL(2,\rr)$ acting on the punctured plane by means of matrix multiplication of vectors. In the case of cocompact subgroups, further properties of such crossed products are deduced from properties of the $C^*$-algebra associated to the horocycle flow on the corresponding compact homogeneous space of $\SL(2,\rr)$.
\end{abstract}

\section{Introduction}
Transformation group $C^*$-algebras represent a tool for the construction of examples of structure and classification theorems for $C^*$-algebras and provide a way to interpret dynamical properties on the $C^*$-algebraic level. Typical examples are the $C^*$-algebras associated to minimal homeomorphisms on infinite compact metric spaces with finite covering dimension (\cite{rieffel-ir,giordano-putnam-skau,toms-winter}), or more generally, free minimal actions of countable residually finite groups with asymptotically finite-dimensional box space on compact metric spaces with finite covering dimension (\cite{swz}). In these cases the structure is that of an $ASH$-algebra and classification is provided by the Elliott invariant. Moving to the non-unital setting, such classification results are still available in the case the $C^*$-algebra is stable and contains projections, assuming a suitable Rokhlin type property for the action. In these situations the resulting transformation group $C^*$-algebra is a stabilized $ASH$-algebra. Examples come from free and minimal actions of the real numbers on compact metric spaces admitting compact transversals (\cite{hsww}), where stability is reminiscent of freeness and the transversal produces a projection in the crossed product. On the other hand, stable simple $\mathcal{Z}$-stable projectionless $C^*$-algebras admit a description of the isomorphism classes of hereditary $C^*$-subalgebras and countable generated Hilbert $C^*$-modules in terms of Cuntz equivalence of positive elements (\cite{zstable_projless}).

Dynamical conditions which ensure stability of a transformation group $C^*$-algebra were given in \cite{green}, where it is proved that $C^*$-algebras arising from actions that are free and wandering on compacts are trivial fields of compact operators. For the case of more general $C^*$-algebras, other characterizations of stability are contained in \cite{rordam-fp,brttw}.\\

The present paper focuses on transformation group $C^*$-algebras associated to the action of discrete subgroups os $\SL(2,\rr)$ on the punctured plane, by means of matrix multiplication of vectors. Ergodic properties of such dynamical systems have been investigated in several places and the duality with the horocycle flow on the corresponding homogeneous spaces for $\SL(2,\rr)$ has been successfully employed in \cite{furstenberg,ledrappier, mau_weiss}. The study of such dynamical systems and their generalizations has a number of interesting applications, as observed in \cite{goro_weiss}, as the quantitative Oppenheim conjecture, quantitative estimates of the denseness of certain projections associated to irreducible lattices and strengthenings of distribution results concerning actions of lattices by automorphisms.

The first part of this work focuses on the study of the distribution of orbits of compact sets on the punctured plane under the action of discrete subgroups of $\SL(2,\rr)$ containing two hyperbolic elements with different axes. Rather than studying the asymptotics of the distribution of such orbits under an increasing family of finite subsets in the lattice, as in \cite{nogueira2002}, \cite{nogueira2010} and \cite{guilloux}, we consider the possibility to find, at every step, an element in the group that \textit{squeezes} enough the image of the compact set under the action of any element in the finite subset. This property of the action resembles the fact that such discrete subgroups of $\SL(2,\rr)$ actually contain an abundance of hyperbolic elements and represents a weaker version of the wandering on compacts assumption considered in \cite{green}. This dynamical condition guarantees the existence of invertible approximants for the elements in the crossed product $C^*$-algebra. By appealing to \cite{rordam-fp}, we show that in the case of actions that are contractive in a suitable sense, this property is enough in order to ensure stability of the crossed product $C^*$-algebra. The "dual" approach is used in the last part to find properties of the crossed product arising from an action of a cocompact subgroup of $\SL(2,\rr)$ on $\rr^2 \backslash \{0\}$ by establishing a $*$-isomorphism between this $C^*$-algebra and the $C^*$-algebra associated to the horocycle flow on the corresponding homogeneous space for $\SL(2,\rr)$.\\

\subsection{Notation}
If $G$ is a locally compact group and $A$ is a $C^*$-algebra, by an action of $G$ on $A$ we mean a continuous group homomorphism from $G$ to the group $\Aut (A)$ of $*$-automorphisms of $A$, endowed with the topology of pointwise convergence. If $X$ is a locally compact Hausdorff space, by an action of $G$ on $X$ we mean a continuous map $G \times X \ra X$ that is associative and such that the identity of the group leaves every point of the space fixed.
If a locally compact group $G$ acts on a locally compact Hausdorff space $X$ by means of an action $\alpha : G \times X \ra X$, we denote by $C_0 (X) \rtimes G$ the associated (full) transformation group $C^*$-algebra, that is the full crossed product $C^*$-algebra relative to the action $\hat{\alpha}_g (f)=f \circ g^{-1}$ for $g \in G$, $f \in C_0 (X)$. Similarly $C_0 (X) \rtimes_r G$ is the reduced transformation group $C^*$-algebra, that is the reduced crossed product relative to the same action.
If $X$ and $Y$ are two Hilbert modules over a $C^*$-algebra, we write $X\Subset Y$ to mean that $X$ is compactly contained in $Y$ in the sense of \cite{cuntz_hm} Section 1.
If $F \subset S$ is an inclusion of sets, we write $F\Subset S$ to mean that $F$ has finite cardinality.
If $X$ is a topological space and $S\subset X$ a subset, we denote by $S^\circ$ its interior.

\section{Weak stable rank $1$}
The concept of stable rank for $C^*$-algebras was introduced by Rieffel in \cite{rieffel} as a noncommutative analogue of the covering dimension of a space and the case of stable rank $1$ is of particular interest (see for example \cite{cuntz_hm} and \cite{open_proj}).
Conditions under which a transformation group $C^*$-algebra has stable rank $1$ have been given in \cite{poon} for actions of the integers; for actions of other groups with finite Rokhlin dimension on compact spaces such conditions can be obtained by combining the results in \cite{hwz}, \cite{szabo} or \cite{swz} and \cite{rordam-sr}, under some other assumptions, as for example, the existence of an invariant measure. If $A$ is a $C^*$-algebra, it is said to have stable rank $1$ if every element in its minimal unitization $\tilde{A}$ can be approximated by invertible elements in $\tilde{A}$. We will consider a more restrictive (non-stable) approximation property, which was used in a crucial way in \cite{brttw}. The following definition was given by Hannes Thiel during a lecture about the Cuntz semigroup in the Winter semester 2016/2017 at the University of M{\"u}nster.
\begin{defn}
\label{defn2.0}
Let $A$ be a $C^*$-algebra. Then $A$ has \textit{weak stable rank $1$}, $\wsr (A)=1$, if $A \subset \overline{GL(\tilde{A})}$.
\end{defn}
Another variation of the concept of stable rank $1$ is the following
\begin{defn}[\cite{zstable_projless} Definition 3.1]
Let $A$ be a $C^*$-algebra. Then $A$ has \textit{almost stable rank $1$}, $\asr (A)=1$, if $\wsr (B) =1$ for every hereditary $C^*$-subalgebra $B \subset A$.
\end{defn}

A $C^*$-algebra $A$ is said to be stable if $A\otimes \mathbb{K} \simeq A$, where $\mathbb{K}$ denotes the $C^*$-algebra of compact operators on a separable Hilbert space. Stable $C^*$-algebras always have weak stable rank $1$ by Lemma 4.3.2 of \cite{brttw} and their multiplier algebra is properly infinite by \cite{rordam-fp} Lemma 3.4. The connection between stability and stable rank in the $\sigma$-unital case was already investigated in \cite{rordam-fp} Proposition 3.5 and Proposition 3.6. For our purpose, we need the following slight variation of the results contained in \cite{rordam-fp}:
\begin{thm}
\label{thm2.0}
Let $A$ be a $\sigma$-unital $C^*$-algebra. The following are equivalent
\begin{itemize}
\item[(i)] $\wsr (A)=1$ and $M(A)$ is properly infinite;
\item[(ii)] $A$ is stable.
\end{itemize}
If $A$ is simple, they are equivalent to
\begin{itemize}
\item[(iii)] $\wsr(A)=1$ and $M(A)$ is infinite.
\end{itemize}
\end{thm}
\proof 
The proof of Lemma 3.2 of \cite{rordam-fp} applies under the hypothesis of weak stable rank $1$, hence if $\wsr (A) =1$ and $M(A)$ is properly infinite, then $A$ is stable by the considerations in the proof of \cite{rordam-fp} Proposition 3.6. As already observed, for any stable $C^*$-algebra $A$, $\wsr(A)=1$ and its multiplier algebra is properly infinite. In the simple case the result follows by an application of Lemma 3.3 of \cite{rordam-fp} and the proof is complete.\\

In order to obtain stability for a transformation group $C^*$-algebra, we introduce a certain dynamical condition and observe that it guarantees weak stable rank $1$; this is the content of the rest of this section. We will deduce infiniteness properties for the multiplier algebra by adapting the results contained in \cite{sth} to the locally compact case in the next section.

\begin{defn}
\label{defn2.1}
Let $G$ be a discrete group acting on a locally compact Hausdorff space $X$. The action is said to be \textit{squeezing} if for every $F\Subset G$ and every $C \subset X$ compact there exists $\gamma \in G$ such that
\[
\gamma g \gamma h C \cap \gamma g C \cap C =\emptyset \nonumber
\]
for all $g,h \in F$.
\end{defn}
Note that Definition \ref{defn2.1} only makes sense for actions on locally compact non-compact spaces, since the space itself is globally fixed by any homeomorphism.

\begin{pro}
\label{prop2.1}
Let $G$ be a discrete group acting on a locally compact Hausdorff space $X$ by means of a squeezing action. Then $\wsr (C_0 (X) \rtimes G) =1$.
\end{pro}
\proof 
Every element in $C_0 (X) \rtimes G$ can be approximated by elements in $C_c (G, C_c (X))$, hence it is enough to prove that any element in $C_c (G, C_c (X))$ is the limit of invertible elements in $(C_0 (X) \rtimes G )^\sim$. Let $F \Subset G$ and $z=\sum_{g \in F} z_g u_g$ be such that $z_g \in C_c (X)$ for every $g \in F$. Define $C := \bigcup_{g \in F} \supp (z_g)$ and let $K \subset X$ be a compact subset such that $C \subsetneq K^\circ$. There exists a continuous function $f : X \ra [0,1]$ such that $\supp (f) \subset K$, $f|_{C} =1$; furthermore, since the action is squeezing, there is a group element $\gamma \in G$ such that 
\[
\gamma g \gamma h K \cap \gamma g K \cap K = \emptyset \nonumber \qquad \mbox{ for all }\quad g,h \in F \cup F^{-1} \cup \{ e\}.
\]
From our choice of $f$, it follows that we can write $z = f z = (f u_{\gamma})(u_{\gamma^{-1}} z)$. 
Computing the third power of $u_{\gamma^{-1}} z$ we obtain
\[
(u_{\gamma^{-1}} z)^3 = \sum_{g,g' , g'' \in F} (z_g \circ \gamma) (z_{g'} \circ (\gamma^{-1} g \gamma^{-1})^{-1}) (z_{g''} \circ (\gamma^{-1} g \gamma^{-1} g' \gamma^{-1})^{-1}) u_{\gamma^{-1} g \gamma^{-1} g' \gamma^{-1} g''}.
 \nonumber
\]
For every $s \in G$ and $\phi \in C_c (X)$ we have $\supp (\phi \circ s^{-1}) = s \supp (\phi)$ and so from our choice of $K$ and $\gamma$, we see that
\[
\begin{split}
&\supp (z_g \circ \gamma) \cap \supp (z_{g'} \circ (\gamma^{-1} g \gamma^{-1})^{-1}) \cap \supp (z_{g''} \circ (\gamma^{-1} g \gamma^{-1} g' \gamma^{-1})^{-1}) \\
&\subset \gamma^{-1} (K \cap  g \gamma^{-1} K \cap g \gamma^{-1} g' \gamma^{-1} K) = \emptyset,
\end{split} \nonumber
\]
since $K \cap  g \gamma^{-1} K \cap g \gamma^{-1} g' \gamma^{-1} K =\emptyset$ if and only if $\gamma(g')^{-1} \gamma g^{-1} K \cap \gamma (g')^{-1} K \cap K = \emptyset$. Hence 
\[
(u_\gamma^{-1} z)^3 = 0  \nonumber
\]
and $u_\gamma z$ is nilpotent. In the same way we obtain
\[
(f u _{\gamma})^3 = f (f\circ \gamma^{-1}) (f\circ \gamma^{-2}) u_{\gamma^{3}}=0 \nonumber
\]
since $\gamma^2 K \cap \gamma K \cap K = \emptyset$. Hence $z$ is a product of nilpotent elements, thus is the limit of invertible elements in $(C_0 (X) \rtimes G)^\sim$ (cfr. \cite{rordam-uhf} 4.1) and the claim follows.

\begin{rem}
\label{oss2.0}
Natural variations of Definition \ref{defn2.1} lead to the same result of Proposition \ref{prop2.1}. The reason why we chose this form is that it fits in the discussion of Section 4.
\end{rem}

\begin{rem}
\label{oss2.1}
Proposition \ref{prop2.1} applies to the reduced crossed product as well.
\end{rem}

\section{Contractive and paradoxical actions}
In the last section we determined a condition on an action of a discrete group that guarantees weak stable rank $1$ for the transformation group $C^*$-algebra. In view of Theorem \ref{thm2.0} this section is devoted to find conditions that guarantee infiniteness properties for the multiplier algebra of the crossed product $C^*$-algebra.\\

If $A$ is any $C^*$-algebra and $G$ a discrete group acting on it, then $A\rtimes G$ is isomorphic to an ideal in $M(A) \rtimes G$, where the action of $G$ on $M(A)$ is the extension of the action on $A$. Then there is a unital $*$-homomorphism $\phi : M(A) \rtimes G \ra M(A \rtimes G)$; if we identify $M(A\rtimes G)$ with the $C^*$-algebra of double centralizers on $A\rtimes G$ and $A\rtimes G$ with its isomorphic image in $M(A) \rtimes G$, $\phi (x) y = xy$ for any $x$ in $M(A) \rtimes G$ and $y$ in $A\rtimes G$. The same results apply to the reduced crossed product as well. This will be the framework for the following considerations.\\

In virtue of the above discussion, all the results we state in the rest of this section concerning full transformation group $C^*$-algebras hold true for the reduced transformation group $C^*$-algebras as well. The same applies to the results contained in the next section, where in order to prove the analogue of Proposition \ref{prop123} for the reduced crossed product, one can use the extension of the surjective $*$-homomorphism from the full crossed product to the reduced crossed product to the multiplier algebras.

The concept of contractive action (see below) was already considered in \cite{sth} page 22 and has to be compared with the more restrictive Definition 2.1 of \cite{delaroche}.
\begin{defn}
\label{defn3.1}
Let $G$ be a discrete group acting on a locally compact Hausdorff space $X$. The action is said to be \textit{contractive} if there exist an open set $U \subset X$ and an element $t \in G$ such that $t \overline{U} \subsetneq U$. In this case $(U,t)$ is called a \textit{contractive pair} and $U$ a \textit{contractive set}.
\end{defn}
The notion of scaling element was introduced in \cite{blackadar-cuntz} and was used to characterize stable algebraically simple $C^*$-algebras.
\begin{defn}[\cite{blackadar-cuntz} Definition 1.1]
\label{defn3.3}
Let $A$ be a $C^*$-algebra and $x$ an element in $A$. $x$ is called a \textit{scaling element} if $x^* x (xx^*) = xx^*$ and $x^* x \neq xx^*$.
\end{defn}
\begin{pro}
\label{prop3.1}
Let $G$ be a discrete group acting on a locally compact Hausdorff space $X$. Consider the following properties:
\begin{itemize}
\item[(i)] The action of $G$ on $X$ is contractive.
\item[(ii)] There exists a scaling elementary tensor in $C_c (X, C_b (X))$.
\end{itemize}
Then $(ii) \Rightarrow (i)$. If $X$ is normal, then $(i) \Rightarrow (ii)$.
\end{pro}
\proof 
$(ii) \Rightarrow (i)$: Let $x=u_t f$ be a scaling elementary tensor in $C_c (G, C_b (X))$ and $U$ the interior of $\supp(f)$. Since $x^* x = |f|^2$ and $xx^* = | f \circ t^{-1} |^2$, the condition $x^* x xx^* = xx^*$ implies $|f| |_{t\overline{U}} =1$; in particular $t\overline{U} \subset U$. Suppose that $t\overline{U} =U$. Then
\[
|f| |_{U^c} =0, \quad |f||_U = |f||_{t\overline{U}} = 1|_{t\overline{U}} \nonumber
\]
and
\[
|f\circ t^{-1} ||_{U^c} = |f\circ t^{-1} ||_{(t\overline{U})^c} =| f\circ t^{-1} ||_{t (U)^c}=0. \nonumber
\]
Since $G$ acts by homeomorphisms, $U$ is a clopen set and $t^{-1} U = t \overline{U}$, which entails
\[
|f\circ t^{-1}||_U =1|_U = 1|_{t \overline{U}}. \nonumber
\]
This would imply $|f|= |f\circ t^{-1}|$ and $x^* x=xx^*$. Hence $t\overline{U} \subsetneq U$.\\
Suppose now that $X$ is normal and let $(U,t)$ be a contractive pair. Take $\xi \in U \backslash (t\overline{U})$. By Urysohn Lemma (normality) there exists a continuous function $f : X \ra [0,1]$ that is $0$ on $U^c$ and $1$ on $\{ \xi \} \cup (t\overline{U})$. The element $x := u_t f \in C_c (G, C_b (X))$ satisfies $x^* x = f^2$, $xx^* = (f\circ t^{-1} )^2$ and $x^* x (xx^*)=xx^*$. Since $\supp (f\circ t^{-1} ) \subsetneq \supp (f)$, we have $x^* x \neq xx^*$, completing the proof.

\begin{cor}
\label{cor3.3.1}
Let $G$ be a group acting on a locally compact normal Hausdorff space by means of a contractive action. Then $M(C_0 (X) \rtimes G)$ is infinite.
\end{cor}
\proof 
Let $x$ be as in Proposition \ref{prop3.1}, we want to show that $\phi(x^*x) \neq \phi (xx^*)$ ($\phi$ is defined at the beginning of the section). For take $\xi \in U \backslash (t\overline{U})$ and let $f \in C_c (X)$ be such that $f(\xi)=1$. Then $(x^*x f)(\xi) \neq 0$ and $(xx^* f)(\xi)=0$ and so $\phi (x^* x) \neq \phi (xx^*)$. As shown in \cite{blackadar-cuntz} Theorem 3.1 the element $\phi(x)+(1-\phi(x^* x))^{1/2}$ is a nontrivial isometry and the claim follows.\\

A variation of the concept of contractive action is the following (see \cite{sth} Lemma 2.3.2) and is a particular case of Definition 2.3.6 of \cite{sth}.
\begin{defn}
\label{defn3.33}
Let $X$ be a locally compact Hausdorff space and $G$ a discrete group acting on it. We say that the action is \textit{paradoxical} if there are positive natural numbers $n$, $m$, group elements $t_1 ,..., t_{n+m}$ and non-empty open sets $U_1 ,..., U_{n+m}$ such that $\bigcup_{i=1}^n U_i = \bigcup_{i=n+1}^{n+m} U_i = X$, $\bigcup_{i=1}^{n+m} t_i (U_i) \subsetneq X$ and $t_i U_i \cap t_j U_j = \emptyset$ for every $i\neq j$.
\end{defn}

Adapting the ideas (and methods) of \cite{sth} Lemma 2.3.7 to the locally compact case, we have the following 
\begin{pro}
\label{prop3.2}
Let $G$ be a discrete group acting on a locally compact normal Hausdorff space $X$. If the action is paradoxical, then $M(C_0 (X) \rtimes G)$ is properly infinite.
\end{pro}
\proof 
Let  $n$, $m$, $t_1 ,..., t_{n+m}$ and $U_1 ,..., U_{n+m}$ be as in Definition \ref{defn3.33}.
Taking unions and relabeling we can suppose $t_i \neq t_j$ for $i\neq j$. Let $F:= \{ t_1 ,..., t_n \}$, $F' := \{ t_{n+1} ,..., t_{n+m}\}$.\\
Since $X$ is normal we can take a partition of unity $\{\phi_t\}_{t \in F}$ subordinated to $\{U_i\}_{i=1}^n$ and a partition of unity $\{ \psi_{s}\}_{s \in F'}$ subordinated to $\{U_i\}_{i=n+1}^{n+m}$. Consider the extension of the action of $G$ to $C_b (X)$ and the associated crossed product $C^*$-algebra $C_b (X) \rtimes G$.\\
Define $x:= \sum_{t \in F} u_t \phi_t^{1/2}$ and $y:= \sum_{t' \in F'} u_{s} \psi_{s}^{1/2}$. Then
\[
x^* x = y^* y = 1. \nonumber
\]
Note now that
\[
x^*y = \sum_{t \in F, s \in F'} \phi_t^{1/2} (\psi_s^{1/2} \circ s^{-1} t ) u_{t^{-1} s} =0 \nonumber
\]
and so $xx^* \perp yy^*$.\\
Let $\phi :C_b (X) \rtimes G \ra M(C_0 (X) \rtimes G)$ be as at the beginning of this section. Take a positive function $f \in C_c(X)$ that takes the value $1$ on a point $\xi \in (\bigcup_{1\leq i \leq n} t_i U_i)^c$. Then
\[
xx^* f = \sum_{t, t' \in F} (\phi_t^{1/2} \circ t^{-1}) (\phi_{t'}^{1/2} \circ t^{-1}) u_{t(t')^{-1}} f \nonumber
\]
entails $0=(xx^* f)(\xi) \neq f (\xi) =1$. Hence $xx^* f \neq f$ and $\phi (xx^*)\neq \phi (1)=1$. The same applies to $yy^*$ and so $1 \in M(C_0 (X) \rtimes G)$ is properly infinite, as claimed.\\

If a discrete group $G$ acts on a locally compact Hausdorff space $X$, the action is said to be \textit{topologically free} if for every $F \Subset G$ the set $\bigcap_{t \in F \backslash \{e\}} \{ x \in X \; | \; tx \neq x \}$ is dense in $X$ (\cite{archbold-spielberg} Definition 1).
Combining Proposition \ref{prop3.2} with the results of Section $2$ we obtain
\begin{thm}
\label{thm3.1}
Let $G$ be a discrete group acting on a locally compact metric space by means of an action that is paradoxical and squeezing. Then $C_0 (X) \rtimes G$ is stable. If the action is topologically free, minimal, squeezing and contractive, then $C_0 (X) \rtimes_r G$ is stable.
\end{thm}
\proof 
Since $X$ is second countable, $C_0 (X) \rtimes G$ is separable, hence $\sigma$-unital. The result follows from Theorem \ref{thm2.0}, Proposition \ref{prop3.2} and Proposition \ref{prop2.1}. If the action is topologically free and minimal, then $C_0 (X)\rtimes_r G$ is simple by \cite{archbold-spielberg}; hence Theorem \ref{thm2.0} applies also in this situation.

\section{The case of discrete subgroups of $\SL(2,\rr)$}
A Fuchsian group $\Gamma$ is a discrete subgroup of $\PSL(2,\rr)$ (\cite{katok} Definition 2.2) and as such it acts on the hyperbolic plane $\mathbb{H}$ and on its boundary $\partial \mathbb{H} =\rr\cup \{\infty\} \simeq \rr \mathbb{P}^1$ by means of M{\"o}bius transformations.
Let $G$ be a discrete subgroup of $\SL(2,\rr)$ acting on $\rr^2 \backslash \{0\}$ by means of matrix multiplication of vectors. The quotient map $\pi : \rr^2 \backslash \{0\} \ra \rr \mathbb{P}^1$ induces an action of $G$ on $\rr\mathbb{P}^1$, which factors through the action of the corresponding Fuchsian group $p (G)$, where $p: \SL(2,\rr) \ra \PSL(2,\rr)$ is the quotient by the normal subgroup $\{-1,+1\}$ of $\SL(2,\rr)$.\\
If $\gamma$ is a hyperbolic element (\cite{katok} 2.1) in $\PSL(2,\rr)$ or $\SL(2,\rr)$ acting on $\mathbb{RP}^1$, we denote by $\gamma^{- (+)}$ its repelling (attracting) fixed point.
For a subset of $\SL(2,\rr)$ or $\PSL(2,\rr)$ consisting of hyperbolic transformations, we say that its elements have different axes if the fixed-point sets for the action of the elements on $\rr \mathbb{P}^1$ are pairwise disjoint. Note that in both $\SL(2,\rr)$ and $\PSL(2,\rr)$, discreteness of a subgroup $G$ implies that whenever two hyperbolic elements in $G$ have a common axis, then both their axes coincide.
\begin{lem}
\label{lem5}
Let $\Gamma$ be a Fuchsian group containing two hyperbolic elements with different axes. Then for every $F\Subset \Gamma$ there exists a hyperbolic element $\gamma \in \Gamma$ such that
\[
g \gamma^+ \neq \gamma^- \qquad \forall g \in F. \nonumber
\]
The same is true if $\Gamma$ is a group generated by a hyperbolic element.
\end{lem}
\proof 
Let $F \Subset \Gamma$. If $\Gamma$ contains two hyperbolic elements with different axes, then it contains infinitely many, hence we can take $\eta$, $\delta$ hyperbolic with different axes and such that the fixed points of $\eta$ are not fixed by any elements in $F$. Suppose that $F$ is such that for every $n \in \nn$ there is a $g \in F$ with $g \eta^n \delta^+ =g(\eta^n \delta \eta^{-n})^+ = (\eta^n \delta \eta^{-n})^- = \eta^n \delta^-$. Then, passing to a subsequence
\[
\exists g \in F \quad \mbox{ s.t. } \quad g\eta^{n_k} \delta^+ = \eta^{n_k} \delta^-. \nonumber
\]
Both $\eta^{n_k} \delta^+$ and $\eta^{n_k} \delta^-$ converge to $\eta^+$ and so $\eta^+$ is fixed by $g$, a contradiction.

\begin{pro}
\label{prop4.1}
Let $G$ be a discrete subgroup of $\SL(2,\rr)$ such that $p(G)$ is a Fuchsian group containing two hyperbolic elements with different axes or a Fuchsian group generated by a hyperbolic element. Then the action of $G$ on $\rr^2 \backslash \{0\}$ is squeezing.
\end{pro}
\proof 
Let $F \Subset G$ and let be given an orthonormal basis $\{e_1 , e_2\}$ for $\rr^2$. By Lemma \ref{lem5} there is a $\gamma \in p (G)$ such that $p(g) \gamma^+ \neq \gamma^-$ for every $g \in G$. Let $h \in G$ be such that $p(h) = \gamma$. Hence $h$ is hyperbolic and is conjugated in $\SL(2,\rr)$ to a diagonal matrix:
\[
h = u^{-1} \Lambda u =u^{-1} \left( \begin{array}{cc}	\lambda	&	0	\\
							0		&	\lambda^{-1}	\end{array}\right)u, \qquad |\lambda| >1. \nonumber
\]
Let $g'$ and $g$ be elements in $G$ and suppose that the upper-left diagonal entry of the matrix $ug'u^{-1}$ vanishes: $(ug'u^{-1})_{1,1} =0$. This means that $\langle e_1, ug'u^{-1}e_1\rangle =0$, or equivalently,  $ug'u^{-1} e_1 \in \rr e_2$; hence, since the image of $u^{-1} e_1$ under the quotient map $\pi : \rr^2 \backslash \{0\} \ra \rr \mathbb{P}^1$ is $\gamma^+$ and the image of $u^{-1} e_2$ under the same map is $\gamma^-$, looking at the action of $p(G)$ on $\rr \mathbb{P}^1$ we obtain $p (g') \gamma^+ = \gamma^-$, contradicting the assumption. Hence $(ug'u^{-1})_{1,1}\neq 0$. Define $g_u := ugu^{-1}$, $g'_u := ug' u^{-1}$ and compute for $n \in \nn$
\[
\Lambda^n g_u= \left(\begin{array}{cc}	\lambda^n (g_{u})_{1,1}		&	\lambda^n (g_u)_{1,2}	\\
								\lambda^{-n} (g_{u})_{2,1}	&	\lambda^{-n}	(g_{u})_{2,2}	\end{array}\right), \nonumber 
\]
and
\[
\Lambda^n g'_u \Lambda^n g_u = \left( \begin{array}{cc}	\lambda^{2n} (g'_u)_{1,1} (g_{u})_{1,1} + (g'_u)_{1,2} (g_u)_{2,1}	&	\lambda^{2n} (g'_u)_{1,1} (g_u)_{1,2} + (g'_u)_{1,2} (g_u)_{2,2} \\
(g'_u)_{2,1} (g_u)_{1,1} + \lambda^{-2n} (g'_u)_{2,2} (g_u)_{2,1} 	&	(g'_u)_{2,1} (g_u)_{1,2}  + \lambda^{-2n} (g'_u)_{2,2} (g_u)_{2,2} \end{array}\right). \nonumber
\]
Let $C \subset \rr^2 \backslash \{ 0 \}$ be a compact subset; take real positive numbers $r_1$ and $r_2$ such that the compact crown $C_{r_1 , r_2} = \{ z  \in \rr^2 | r_1 \leq \| z \| \leq r_2\}$ contains $uC$. We want to show that there exists $n >0$ such that
\[
\Lambda^n g'_u \Lambda^n g_u C_{r_1 , r_2} \cap \Lambda^n g'_u C_{r_1 , r_2} \cap C_{r_1,r_2} = \emptyset, \qquad \forall g,g' \in F. \nonumber
\]
Let $(x,y)^t \in \rr^2$ be such that $\Lambda^n g'_u \Lambda^n g_u (x,y)^t$ belongs to $uC$. In particular, this entails
\[
\| \Lambda^n g'_u \Lambda^n g_u (x,y)^t \| \leq r_2 \nonumber
\]
and taking the first coordinate:
\[
| \lambda^{2n} (g'_u)_{1,1} [(g_u)_{1,1} x + (g_u)_{1,2} y] + [ (g'_u)_{1,2} (g_u)_{2,2} x + (g'_u)_{1,2} (g_u)_{2,2} y ]| \leq r_2. \nonumber
\]
Hence
\begin{equation}
\label{eq4.1}
|(g_u)_{1,1} x + (g_u)_{1,2} y| \leq \frac{r_2 + | (g'_u)_{1,2} (g_u)_{2,2} x + (g'_u)_{1,2} (g_u)_{2,2} y |}{\lambda^{2n} | (g'_u)_{1,1}|}
\end{equation}
for every $(x,y)^t \in (\Lambda^n g'_u \Lambda^n g_u)^{-1} C$. Furthermore, if $(x,y)^t \in \rr^2$ is such that $\Lambda^n g_u (x,y)^t $ belongs to $C_{r_1,r_2}$, then
\begin{equation}
\label{eq4.2}
\begin{split}
r_1^2 &\leq [\lambda^n (g_u)_{1,1} x + \lambda^n (g_u)_{1,2} y]^2 + [\lambda^{-n} (g_u)_{2,1} x + \lambda^{-n} (g_u)_{2,2} y ]^2\\
& = \lambda^{2n} [(g_u)_{1,1} x + (g_u)_{1,2} y]^2 + \lambda^{-2n}[ (g_u)_{2,1} x + (g_u)_{2,2} y ]^2.
\end{split}
\end{equation}
Combining (\ref{eq4.1}) and (\ref{eq4.2}) we obtain
\begin{equation}
\label{eq4.3}
r_1^2 \leq  \frac{[r_2 + | (g'_u)_{1,2} (g_u)_{2,2} x + (g'_u)_{1,2} (g_u)_{2,2} y |]^2}{\lambda^{2n} | (g'_u)_{1,1}|^2} + \lambda^{-2n} [(g_u)_{2,1} x + (g_u)_{2,2} y]^2 
\end{equation}
for every $(x,y)^t \in (\Lambda^n g'_u \Lambda^n g_u)^{-1} C_{r_1,r_2} \cap (\Lambda^n g_u)^{-1} C_{r_1,r_2}$.
If $(x,y)^t$ belongs to $C_{r_1,r_2}$, then there is a constant $M >0$ such that
\[
\frac{| (g'_u)_{1,2} (g_u)_{2,2} x + (g'_u)_{1,2} (g_u)_{2,2} y |]^2}{ | (g'_u)_{1,1}|^2} +  [(g_u)_{2,1} x + (g_u)_{2,2} y]^2 \leq M \nonumber
\]
and this constant does not depend on the choice of $g$, $g'$ in $F$. So, by (\ref{eq4.3}), for $n$ large enough
\[
(\Lambda^n g'_u \Lambda^n g_u)^{-1} C_{r_1,r_2}  \cap (\Lambda^n g_u)^{-1} C_{r_1,r_2} \cap C_{r_1,r_2} = \emptyset, \nonumber
\]
which entails
\[
C_{r_1,r_2} \cap \Lambda^n g'_u C_{r_1,r_2} \cap \Lambda^n g'_u \Lambda^n g_u C_{r_1,r_2}=\emptyset \nonumber
\]
and so
\[
u^{-1}C_{r_1,r_2} \cap h^n g' u^{-1}C_{r_1,r_2} \cap h^n g' h^n g u^{-1}C_{r_1,r_2}=\emptyset. \nonumber
\]
The result follows since $C \subset u^{-1} C_{r_1 , r_2}$. \\

Hence we have determined a class of discrete subgroups of $\SL(2,\rr)$ whose action on $\rr^2 \backslash \{0\}$ is squeezing. Conditions under which this action is contractive or paradoxical are the content of the following

\begin{pro}
\label{prop4.2}
Let $G$ be a discrete subgroup of $\SL(2,\rr)$ acting on $\rr^2 \backslash \{ 0 \}$ by means of matrix multiplication of vectors. If $G$ contains a hyperbolic element, then the action is contractive. If $G$ contains at least two hyperbolic elements with different axes, then the action is paradoxical.
\end{pro}
\proof
Suppose that $G$ contains a hyperbolic element, then the same is true for its image under the quotient map $p : \SL(2,\rr) \ra \PSL(2,\rr)$. Since the action of $\Gamma = p (G)$ on $\rr\mathbb{P}^1$ is by homeomorphisms and every hyperbolic element in $\Gamma$ is conjugated in $\PSL(2,\rr)$ to a Moebius transformation of the form $z \mapsto \lambda^2 z$ for some $\lambda >1$, it follows that the action of $\Gamma$ on $\rr\mathbb{P}^1$ is contractive. Hence there are $U \subset \rr\mathbb{P}^1$ and $\gamma \in \Gamma$ such that
\begin{equation}
\gamma \overline{U} \subsetneq U. \nonumber
\end{equation}
In the case $G$ contains at least two hyperbolic elements with different axes, then the same is true for $\Gamma$ and as is well known, in this case $\Gamma$ contains a countable subset of hyperbolic elements with different axes. In order to see this, let $\gamma$, $\eta$ be hyperbolic elements in $\Gamma$ with different axes; then the elements in the sequence $\{ \eta^n \gamma \eta^{-n}\}_{n \in \nn}$ are hyperbolic transformations with different axes. In particular, for every $n, m \geq 2$ natural numbers there are group elements $\gamma_1 ,..., \gamma_{n+m}$ and contractive open sets $U_1 ,..., U_{n+m}$, where for each $i=1,...,n+m$ $U_i$ contains the attracting fixed point $\gamma_i^+$ of $\gamma_i$, such that
\begin{equation}
\label{eq2}
\bigcup_{i=1}^n U_i = \bigcup_{j=n+1}^{n+m} U_j = \rr\mathbb{P}^1,
\end{equation}
\begin{equation}
\label{eq3}
\gamma_i U_i \cap \gamma_j U_j =\emptyset \qquad \forall i\neq j.
\end{equation}
Hence, we just need to observe that the same holds after replacing the sets $U_i$ with $\pi^{-1} (U_i)$ and the elements $\gamma_i$ with some representatives in $G$. Equation (\ref{eq2}) automatically holds for the sets $\pi^{-1} (U_i) \subset \mathbb{R}^2 \backslash \{0\}$. Choose a representative $g_i \in G$ for every $\gamma_i \in \Gamma$; since the action of $G$ on $\rr\mathbb{P}^1$ factors through the action of $\Gamma$, equation (\ref{eq3}) can be replaced by
\[
g_i U_i \cap g_j U_j =\emptyset \qquad \forall i\neq j. \nonumber
\]
By equivariance of the quotient map $\pi : \rr^2 \backslash \{0\} \ra \rr \mathbb{P}^1$ it follows that
\[
g_i (\pi^{-1} (U_i)) \cap g_j (\pi^{-1} (U_j)) = \emptyset \qquad \forall i\neq j. \nonumber
\]
We are left to check that the inverse image of a contractive open set is again a contractive open set. Since the map $\rr^2 \backslash \{0\} \ra \rr\mathbb{P}^1$ is a quotient by a group action (the group is $\rr^\times$), it is open and so the inverse image of the closure of a set is the closure of the inverse image of the same set; hence, if $(U, g)$ is a contractive pair with $U \subset \rr\mathbb{P}^1$ and $g \in G$, then
\[
g \overline{( \pi^{-1} (U))} = g \pi^{-1} (\overline{U}) = \pi^{-1} (g \overline{U}) \subsetneq \pi^{-1} (U). \nonumber
\]
The proof is complete.

\begin{cor}
\label{cor4.1}
Let $G$ be a discrete subgroup of $\SL(2,\rr)$ such that $p (G) \subset \PSL(2,\rr)$ is a Fuchsian group containing two hyperbolic elements with different axes. The transformation group $C^*$-algebra $C_0 (\rr^2 \backslash \{0\})\rtimes G$ is stable.\\
If $p(G)$ is generated by a hyperbolic transformation, then $\wsr (C_0 (\rr^2 \backslash \{0\})\rtimes G) =1$ and $M(C_0 (\rr^2 \backslash \{0\})\rtimes G)$ is infinite.
\end{cor}
\proof 
Follows from Proposition \ref{prop4.2}, Proposition \ref{prop4.1}, Theorem \ref{thm3.1} and Corollary \ref{cor3.3.1}.\\

Corollary \ref{cor4.1} applies to the case of discrete subgroups of $\SL(2,\rr)$ associated to Fuchsian groups of the first kind (\cite{katok} 4.5), hence in particular the cocompact ones. Non-lattice subgroups to which Corollary \ref{cor4.1} applies are considered in \cite{semenova}.\\

In Proposition \ref{prop4.2} we deduced paradoxicality for the action of a discrete subgroup $G$ of $\SL(2,\rr)$ on $\rr^2 \backslash \{0\}$ from paradoxicality of the action of the corresponding Fuchsian group on $\rr \mathbb{P}^1$ and concluded from this fact that the multiplier algebra of $C_0 (\rr^2 \backslash \{0\})\rtimes G$ is properly infinite. 
It follows from \cite{glasner} Example VII.3.6 that if $\Gamma$ is a Fuchsian group of the first kind, then its action on $\rr \mathbb{P}^1$ is extremely proximal (see \cite{glasner} page 96 for the definition) and this property represents a stronger form of paradoxicality, hence stronger infiniteness properties for the multiplier algebra of the transformation group $C^*$-algebra are expected in this case. Note that in \cite{boundary}  an extremely proximal action is called a strong boundary action.\\ The next Proposition is a consequence of the results contained in \cite{boundary} and \cite{kra}.

\begin{lem}
\label{lem4.4}
Let $G$ and $H$ be locally compact groups and $A$, $B$ be $C^*$-algebras. Suppose $G$ acts on $A$ and $H$ acts on $B$ and that there is an equivariant involutive homomorphism $\phi : C_c (G,A) \ra C_c (H,B)$ which is continuous for the $L^1$-norms. Then there is a $*$-homomorphism $\hat{\phi} : A\rtimes G \ra B\rtimes H$.
\end{lem}
\proof 
If $\rho: L^1 (H,B) \ra \mathfrak{H}$ is a nondegenerate $L^1$-continuous involutive representation of $L^1 (H, B)$, then the composition $\rho \circ \phi : C_c (G, A) \ra \mathfrak{H}$ is $L^1$-continuous as well. Hence $\| \phi (f) \| \leq \| f \|$ for every $f \in C_c (G,A)$ by \cite{williams} Corollary 2.46, as claimed.

\begin{pro}
\label{prop123}
Let $G$ be a discrete subgroup of $\SL(2,\rr)$ such that $p(G) \subset \PSL(2,\rr)$ is a finitely generated Fuchsian group of the first kind not containing elements of order $2$. Then $M(C_0 (\rr^2 \backslash \{0\})\rtimes G)$ contains a Kirchberg algebra in the UCT class as a unital $C^*$-subalgebra.
\end{pro}
\proof 
The quotient map $p : \rr^2 \backslash \{0\} \ra \rr \mathbb{P}^1$ is surjective and equivariant with respect to the action of $G$, hence it induces a unital $*$-homomorphism $ C(\mathbb{R} \mathbb{P}^1) \rtimes G \ra C_b (\rr^2 \backslash \{0\} )\rtimes G$ which can be composed with the unital $*$-homomorphism $C_b (\rr^2 \backslash \{0\})\rtimes G \ra M(C_0 (\rr^2 \backslash \{0\}) \rtimes G)$ introduced at the beginning of Section 3 in order to obtain a unital $*$-homomorphism $\phi : C(\mathbb{R} \mathbb{P}^1)\rtimes G \ra M(C_0 (\rr^2 \backslash \{0\}) \rtimes G)$. By \cite{kra}, finitely generated Fuchsian groups of the first kind not admitting elements of order $2$ lift to $\SL(2,\rr)$. Denote by $\kappa : \Gamma \ra  \kappa (\Gamma) \subset G$ a lift. Since the action of $G$ on $\mathbb{R}\mathbb{P}^1$ factors through the action of $\Gamma$, the map
\[
\psi_c : C_c (\Gamma , C(\mathbb{R} \mathbb{P}^1)) \ra C_c (\kappa (\Gamma) , C(\mathbb{R}\mathbb{P}^1)) \nonumber
\]
\[
f \mapsto f \circ \kappa^{-1} \nonumber
\]
is an involutive homomorphism and it preserves the $L^1$-norm, as well as the inclusion $C_c (\kappa (\Gamma), C(\rr \mathbb{P}^1)) \ra C_c (G, C(\rr \mathbb{P}^1))$. By Lemma \ref{lem4.4} there is a (unital) $*$-homomorphism $\psi : C(\mathbb{R}\mathbb{P}^1)\rtimes \Gamma \ra C(\mathbb{R}\mathbb{P}^1) \rtimes G$. Hence $\phi \circ \psi: C(\rr \mathbb{P}^1 ) \rtimes \Gamma \ra M(C_0 (\rr^2 \backslash \{0\})\rtimes G)$ is a unital $*$-homomorphism. By \cite{boundary} Theorem 5 the $C^*$-algebra $C(\rr \mathbb{P}^1 ) \rtimes \Gamma$ is a unital Kirchberg algebra in the UCT class, hence $\psi \circ \phi$ is injective and the result follows.

\section{Cocompact subgroups of $\SL(2,\rr)$}

Consider the one-parameter subgroup of $\SL(2,\rr)$
\[
N:= \{ n(t) \in \SL(2,\rr) \; | \; n(t)=\left( \begin{array}{cc}	1	&	t	\\
										0	&	1	\end{array}\right), \quad t \in \rr \}. \nonumber
\]
Given a discrete subgroup $G$ of $\SL(2,\rr)$, one can define a flow on the corresponding homogeneous space $G \backslash \SL(2,\rr)$ by $Gg \mapsto Gg n(-t)$; this is called the \textit{horocycle flow} (\cite{ew} 11.3.1). The stabilizer of the point $(1,0)^t$ in $\rr^2 \backslash \{0\}$ for the action of $\SL(2,\rr)$ is $N$ and so the quotient $\SL(2,\rr)/N$, endowed with the action of $\SL(2,\rr)$ given by left multiplication, is isomorphic, as a dynamical system, to $\rr^2 \backslash \{0\}$. The interplay between the action of $G$ on $\mathbb{R}^2 \backslash \{0\}$ and the horocycle flow on $G \backslash \SL(2,\mathbb{R})$ is employed in the following

\begin{pro}
\label{prop2.1}
Let $G$ be a discrete cocompact subgroup of $\SL(2,\rr)$. The transformation group $C^*$-algebra $C_0 (\rr^2 \backslash \{0\}) \rtimes G$ is simple, separable, stable, $\mathcal{Z}$-stable, with a unique lower semicontinuous $2$-quasitrace and it has almost stable rank $1$. In particular it satisfies the hypothesis of \cite{io} Theorem 3.5.
\end{pro}
\proof 
Since $G$ is countable and $\mathbb{R}^2 \backslash \{0\}$ is a locally compact second countable Hausdorff space, $C_0 (\mathbb{R}^2 \backslash \{0\}) \rtimes G$ is separable.\\
As already observed in the discussion after Corollary \ref{cor4.1}, $C_0 (\mathbb{R}^2 \backslash \{0\}) \rtimes G$ is stable in the case $G$ is cocompact. Since the action of $G$ on $\rr^2 \backslash \{0\}$ is free and minimal (\cite{ergtopdyn} Theorem IV.1.9), $C_0 (\rr^2 \backslash \{0\}) \rtimes G$ is simple (\cite{archbold-spielberg}). 
By simplicity, the non-trivial lower semicontinuous traces on $C_0 (\mathbb{R}^2 \backslash \{0\}) \rtimes G$ are semifinite (\cite{dixmier} 6.1.3) and so, in virtue of \cite{green2} Proposition 25 and Proposition 26, the restriction map sets up a bijection with the lower semicontinuous semifinite $G$-invariant traces on $C_0 (\rr^2 \backslash \{0\})$. Every such trace is uniquely given by integration against a $G$-invariant Radon measure. By Furstenberg Theorem (\cite{furstenberg}) there is exactly one such non-trivial measure. Hence $C_0 (\mathbb{R}^2 \backslash \{0\}) \rtimes G$ admits a unique non-trivial lower semicontinuous trace. Since the action of $G$ on $\rr^2 \backslash \{0\}$ is amenable, $C_0 (\rr^2 \backslash \{0\}) \rtimes G$ is nuclear (\cite{delaroche2} Theorem 3.4). By exactness it admits a unique non-trivial lower semicontinuous $2$-quasitrace (\cite{kirchberg}).\\
In virtue of Corollary 9.1 and Corollary 6.7 of \cite{hsww} the $C^*$-algebra $C(G \backslash \SL(2,\rr)) \rtimes N$ is stable; hence it follows from Green's imprimitivity Theorem (\cite{williams} Corollary 4.11) that $C_0 (\mathbb{R}^2 \backslash \{0\}) \rtimes G \simeq C(G \backslash \SL(2,\rr)) \rtimes N$. By \cite{hsww} Corollary 9.1 and Theorem 3.5, $C(G \backslash \SL(2,\rr) )\rtimes N$ has finite nuclear dimension; hence, \cite{tikusis} Corollary 8.7 entails $\mathcal{Z}$-stability.\\
As observed in \cite{noncommgeom} page 129, $C(G\backslash \SL(2,\rr)) \rtimes N$ is projectionless; hence, \cite{zstable_projless} Corollary 3.2 applies and $\asr (C(G\backslash \SL(2,\rr)) \rtimes N) =1$.\\
The result follows since the Cuntz semigroup of a stable $\mathcal{Z}$-stable $C^*$-algebra is almost unperforated (\cite{rordam-sr} Theorem 4.5).\\

\begin{rem}
The stability of $C_0 (\rr^2 \backslash \{0\}) G$ in Proposition \ref{prop2.1} can also be established directly from that of $C(G \backslash \SL(2,\rr)) \rtimes N$. In fact, the rest of the proof shows that $C(G \backslash \SL(2,\rr)) \rtimes N$ satisfies the hypothesis of \cite{io} Theorem 3.5. Since $C_0 (\rr^2 \backslash \{0\}) \rtimes G$ is a hereditary $C^*$-subalgebra of $C(G \backslash \SL(2,\rr)) \rtimes N$, it is then enough to prove that the non-trivial lower semicontinuous trace on $C_0 (\mathbb{R}^2 \backslash \{0\}) \rtimes G$ is unbounded. But this follows since it is induced by the Lebesgue measure.
\end{rem}

As a consequence we obtain the following properties for the $C^*$-algebra associated to the action of a cocompact discrete subgroup of $\SL(2,\rr)$ on $\rr^2 \backslash\{0\}$

\begin{cor}
\label{horo_1}
Let $G$ be a cocompact discrete subgroup of $\SL(2,\rr)$, $\tau$ the lower semicontinuous trace associated to the Lebesgue measure $\mu_L$ on $\rr^2 \backslash \{0\}$ and $d_\tau$ the corresponding functional on the Cuntz semigroup $Cu (C_0 (\rr^2 \backslash \{0\})\rtimes G)$. Then
\[
\Ped (C_0 (\rr^2 \backslash \{0\}) \rtimes G)=\{ x \in C(\rr^2 \backslash \{0\}) \rtimes G\; : \; d_\tau ([|x|]) < \infty\}. \nonumber
\]
Every hereditary $C^*$-subalgebra of $C(\rr^2 \backslash \{0\}) \rtimes G$ is either algebraically simple or isomorphic to $C(\rr^2 \backslash \{0\}) \rtimes G$.
\end{cor}

\begin{cor}
\label{horo_3}
Let $G$ be a cocompact discrete subgroup of $\SL(2,\rr)$. Every countably generated right Hilbert module for $C_0 (\rr^2 \backslash \{0\}) \rtimes G$ is isomorphic to one of the form
\[
\overline{f \cdot (C_0 (\rr^2 \backslash \{0\}) \rtimes G)}, \qquad f \in C_0 (\rr^2 \backslash \{0\}). \nonumber
\]
For two such Hilbert modules we have
\[
\overline{f \cdot (C_0 (\rr^2 \backslash \{0\}) \rtimes G)} \simeq \overline{g \cdot (C_0 (\rr^2 \backslash \{0\})\rtimes G)} \qquad \Leftrightarrow \qquad \mu_L (\supp (f)) = \mu_L (\supp (g)) \nonumber
\]
and there exists a Hilbert module $E$ such that
\[
\overline{f \cdot (C_0 (\rr^2 \backslash \{0\})\rtimes G)} \simeq E \Subset \overline{g \cdot (C_0 (\rr^2 \backslash \{0\})\rtimes G)} \nonumber
\]
if and only if $\mu_L (\supp (f)) < \mu_L (\supp (g))$.
\end{cor}
\proof 
The Cuntz semigroup of the $C^*$-algebra $C_0 (\rr^2 \backslash \{0\})\rtimes G$ is stably finite by \cite{cuntz_t} Proposition 5.2.10, hence by \cite{cuntz_t} Proposition 5.3.16 it does not contain compact elements, since this $C^*$-algebra is projectionless. Hence the countably generated Hilbert modules correspond to soft elements and Cuntz equivalence of soft elements is implemented by the unique (up to scalar multiples) nontrivial functional associated to the unique (up to scalar multiples) lower semicontinuous trace $\tau$. It follows that all the possible values in the range of the dimension function are obtained by Cuntz equivalence classes of elements in $C_0 (\mathbb{R}^2 \backslash \{0\})$ since, for every $f \in C_0 (\mathbb{R}^2 \backslash \{0\})$, we have $d_\tau (f) = \mu_L (\supp (f))$. The result follows from Theorem 3.5 of \cite{io}.

\section{Final remarks}
It follows from the results in the last section that if $G$ is a cocompact discrete subgroup of $\SL(2,\rr)$, the Cuntz classes of elements in the transformation group $C^*$-algebra $C_0 (\rr^2 \backslash\{0\})\rtimes G$ are generated by continuous functions on the plane. It might be possible that this property can be derived from the dynamics. 

It can be shown that if we restrict to discrete subgroups of $\SL(2,\rr)$ which are the inverse images under the quotient map $p: \SL(2,\rr) \ra \PSL(2,\rr)$ of fundamental groups of hyperbolic Riemann surfaces, the construction of the $C^*$-algebra associated to the horocycle flow on the corresponding homogeneous space of $\SL(2,\rr)$ induces a functor from a category whose objects are hyperbolic Riemann surfaces and the morphisms are finite sheeted holomorphic coverings to the usual category of $C^*$-algebras; this suggests that it might be possible do detect the holomorphic structure at the $C^*$-algebraic level. Observe that, if $\mathcal{M}_g$ is a compact Riemann surface of genus $g$, after identifying $p^{-1} (\pi_1 (\mathcal{M}_g)) \backslash \SL(2,\rr)$ with the unit tangent bundle $T_1 (\mathcal{M}_g)$, the Thom-Connes isomorphism (\cite{thom-connes}) gives a way to compute the $K$-theory of $C_0 (\rr^2 \backslash \{0\})\rtimes p^{-1}(\pi_1 (\mathcal{M}_g)) \simeq C(T_1 (\mathcal{M}_g))\rtimes \rr$ and it reads
\[
K_0 (C_0 (\rr^2 \backslash \{0\})\rtimes p^{-1}(\pi_1 (\mathcal{M}_g))) = \mathbb{Z}^{2g+1}, \nonumber
\]
\[
K_1 (C_0 (\rr^2 \backslash \{0\})\rtimes p^{-1}(\pi_1 (\mathcal{M}_g)))= \mathbb{Z}^{2g+1} \oplus \mathbb{Z} /(2g-2) . \nonumber
\]
Both the order and the scale in $K_0$ are trivial since $C_0 (\rr^2 \backslash \{0\})\rtimes p^{-1}(\pi_1 (\mathcal{M}_g))$ is projectionless and stable.\\
Furthermore, by \cite{thom-connes} Corollary 2 the range of the pairing between $K_0$ and the unique trace is determined by the Ruelle-Sullivan current associated to this flow (see \cite{noncommgeom} 5-$\alpha$), which is trivial by \cite{paternain}. Thus the Elliott invariant contains information only about the genus, or equivalently, the homeomorphic class of the Riemann surface. In particular, if the Elliott conjecture holds true for this class of $C^*$-algebras and if it is possible to detect the holomorphic structure at the level of the $C^*$-algebras, this should correspond to something finer than the $C^*$-algebraic structure. 
\section{Aknowledgements}
The author thanks Prof. Wilhelm Winter for the hospitality at the Westf\"alische Wilhelms-Universit\"at of M\"unster and Prof. Roberto Longo for the hospitality at the Universit\`a degli Studi di Roma Tor Vergata for the period of this research. Many thanks go to Prof. Ludwik D\k abrowski who carefully read and gave his important feedback on the parts of this paper that are contained in the author's PhD thesis. The author also thanks the anonymous referee for the valuable comments on a previous version of the manuscript which led to an improved exposition. This research is partially supported by INdAM.


\end{document}